
\magnification=\magstep1
\bigskip
\centerline {\bf Automorphisms of crepant resolutions}
\centerline {\bf  for  quotient spaces}
 \bigskip \bigskip
\centerline {{\bf by A.Libgober}%
\footnote{${}^*$}{
the author gratefully
acknowledges  support and hospitality of IHES during preparation of this
paper.}}
\bigskip \bigskip
\centerline {Department of Mathematics}
\centerline {University of Illinois at Chicago}
\centerline {851 S.Morgan Str. Chicago, Ill. 60607}
\centerline {e-mail: libgober@math.uic.edu}
\bigskip \bigskip 1. {\bf Introduction}
\bigskip Let $X$ be a Kahler manifold and $G$ be a finite group
of biholomorphic volume preserving automorphisms of $X$. Let
 $\pi: \widetilde {X/G} \rightarrow X/G $ be a crepant resolution of
the orbit space $X/G$ (i.e. the pullback
of the dualizing sheaf of  $X/G$ is the canonical bundle of
$\widetilde {X/G}$). The euler characteristic of
$ \widetilde {X/G} $ can be found as  follows: (cf. [R],[BD],[HH]):
 $$e(\widetilde {X/G})={1 \over {\vert G \vert}}
 \Sigma_{hg=gh} e(X^g \cap X^h)=
\Sigma_{[g]} e(X^g/C(g)) \eqno (1.0)$$
 Here $X^g$, for $g \in G$, is the fixed point set of an automorphism
and $C(g)$ is the centralizer. The summation in (1.0) is over all
conjugacy classes $[g]$ of elements $g \in G$.
 This formula plays a key role in  validating
that  pairs of Calabi Yau manifolds constitute  mirror pairs (cf. [COGP]).
\par In this paper we shall propose a formula designed for the study of
manifolds for which the pair $(X,G)$ admits a symmetry $h$. This formula
 allows one to calculate the Lefschetz number of an automorphism
acting on a crepant resolution of the quotient (for $h=1$ this
formula becomes (1.0)). We show that it can be derived from a certain
equivariant form of McKay correspondence between a data from a crepant
resolution of quotient singularity ${\bf C}^n /G$, where $G$ admits
certain  automorphism $h$ and the action of $h$ on the conjugacy classes
of $G$. This version of McKay correspondence is valid in dimension $n=2$
 (cf. section 4)
and follows from the existence of certain triangulations for abelian $G$
for $n \ge 3$ (such triangulation can always be constructed in dimension $n=3$
(cf. section 5). In particular this can be used to compare the
Lefschetz numbers of certain involutions on some Calabi Yau threefolds
(considered in  [COGP] and [LT]) and their mirrors (section 3).
Note that detailed investigation of connections between
McKay correspondence and crepant resolutions was carried out in [BD].
\bigskip
More precisely we have the following (here $L(h,Z)$ denotes
the Lefschetz number $\Sigma_i (-1)^i tr(h,H^i (Z))$ of a transformation
$h$ acting on a topological space $Z$):
\bigskip {\bf Theorem 1.}
Let $X$ be a Kahler manifold on which a finite group $G$ of volume preserving
automorphisms
acts holomorphically and let $h$ be a biholomorphic automorphism
of $X$ such that the group $H$
generated by $G$ and $h$ is finite and contains $G$ as a normal
subgroup. Let $\widetilde {X/G}$ be a crepant resolution of $X/G$.
Let us consider the class ${\cal C}(h)$ of elements in $G$ with the
property:
\bigskip (*) In  any $h$-invariant
stabilizer $S_i$ of any point of $X^g$,  the conjugacy class of
$g$ is $h$-invariant.
\bigskip Let us assume that:
\par 1. For $g \in {\cal C}(h)$ $h$ leaves invariant the set
 of stabilizers containing $g$.
\par 2. $h$ normalizes the
centralizer of any element from ${\cal C} (h)$.
\par Then
$$L(h,\widetilde {X/G})=\Sigma_{[g],g \in {\cal C}(h)} L(h,X^g/C(g))
 \eqno (1.1)$$
\bigskip {\bf Remarks}.
 1. Let ${\cal S}$ be the collection of subgroups of $G$ which are
stabilizers of points of $X$. We have: $X^g=
\bigcup_{S_i \in {\cal S}, g \in S_i} X^{S_i}$. Therefore condition 1
implies that if $g \in C(h)$ then $h$ preserves $X^g$.
Also $h$ acts on the quotient $X^g/C(g)$ as a consequence of
condition 2.
\par 2. Throughout the paper we use the following additivity of the
Lefschetz numbers: if $Z=Z_1 \cup Z_2, Z_1 \cap Z_2 =\empty$ is
an $h$ equivariant decomposition then $L(h,Z)=L(h,Z_1) + L(h,Z_2)$ (cf. [Br]).
\par 3. In what follows by stratification we mean just a partition into a
union of disjoint sets.
\par 4. The proof is based on the following statement (an equivariant
form of weak McKay correspondence):
\bigskip (**)
 Let $H$ be a finite subgroup of $GL_n({\bf C})$, $G$ a normal subgroup of
 $H$ which belongs to $SL_n({\bf C})$ and such that $H/G$ is a cyclic group
with a generator $h$. (In particular $h$ induces automorphisms of ${\bf C}^n/G$
and $G$). Then the Lefschetz number of $h$ acting on an $h$-equivariant
crepant resolution of ${\bf C}^n/G$ is equal to the number of $h$-invariant
conjugacy classes in $G$.
\bigskip
 This will be proven below in the case $n=2$ and $n=3$ and $G$ abelian.
\bigskip \bigskip {\bf 2.Proof of Theorem 1}.
\bigskip Let $S_i \in {\cal S}$ ($i=1,..,card {\cal S}$) be collection
of subgroups of $G$ which are stabilizers of points of $X$.
The manifold $X$ admits the stratification $X=\cup X^{S_i}$  such that
each stratum consists of points having the same stabilizer $S_i$.
Let  $X^{[S_i]}$ be the $G$-orbit of the strata
 $X^{S_i}$'s i.e. $X^{[S_i]}$ is the unions of $X^{S_i}$ with $G$-conjugate
$S_i$'s ($[S_i] \in {\cal S}/G$ where $G$ acts on ${\cal S}$ by conjugation).
 The sets $X^{[S_i]}$ form a stratification of
 $X$ with $G$-invariant strata
 and the $G$-quotients  $X^{[S_i]}/G$ of these strata
provide a stratification of
$X/G$.  Let $\pi^{-1}(X^{[S_i]}/G) \subset \widetilde {X/G}$
be the preimages of strata $X^{[S_i]}/G$ in a chosen crepant
resolution $\pi: \widetilde {X/G} \rightarrow X/G$
of $\widetilde {X/G}$ and $\widetilde {X^{[S_i]}/G(h)}$ be
the $h$-orbit of $\pi^{-1}(X^{[S_i]}/G)$.
 The automorphism  $h$ acts on $\widetilde {X^{S_i}/G(h)}$ and we  have:
 $$L(h, \widetilde {X/G}))=\Sigma_{\widetilde {X^{[S_i]}/G(h)} }
 L(h,\widetilde {X^{[S_i]}/G(h)}) \eqno (2.0)$$
 where the summation is over all $h$ orbits of preimages of the
strata $X^{[S_i]}/G$.
\bigskip Next we have:
 $$\Sigma_{\widetilde {X^{[S_i]}/G(h)} }
 L(h,\widetilde {X^{[S_i]}/G(h)})=
 \Sigma_{[S_i] \in {\cal S}, h([S_i])=[S_i])}
 L(h,\widetilde {X^{[S_i]}/G})=$$
$$=\Sigma _{[S_i] \in {\cal S},h([S_i]=[S_i]}
 L(h, {X^{S_i}}/G) \cdot con(h,S_i) \eqno (2.1)$$
 where $con(h,S_i)$ is the number of $h$ invariant conjugacy classes in
 a subgroup $S_i$ and the summation on the last two sums is over
 $h$ invariant conjugacy classes subgroups $[S_i]$.
 (Note that $h([S_i])=[S_i]$ implies that $h$ acts on corresponding
 stratum $X^{[S_i]}$).
 The first equality  takes place since the Lefschetz number of $h$
acting on an orbit of $h$ of a conjugacy class for which
 $h([S_i]) \ne [S_i]$ is zero due to the absence of fixed points.
 The second equality is a consequence of the equivariant McKay correspondence
 as stated in Remark 4.  Note that this remark also implies that the number
 $con(h,S)$ is the same for the conjugate
subgroups (since both singularities $X/S$ and $X/gSg^{-1}$ are
$h$-equivariantly  equivalent to the same singularity of $X/G$).
 The last term can be rewritten as:
  $$\Sigma_{[S_i],S_i \in {\cal S}, h([S_i])=[S_i]}
   L(h, X^{[S_i]}) \cdot {{\vert S_i \vert } \over {\vert G
 \vert }} \cdot con (h,S_i) \eqno (2.2)$$
 where the summation is over all $h$-invariant conjugacy classes
 of stabilizers.  The latter can be replaced by the
sum: $$\Sigma L(h,X^{[S_i]}) \cdot {{\vert S_i \vert }
  \over {\vert G \vert}} \eqno (2.3)$$
 where the summation is over all pairs $([g],[S_i])$ and where $g$ is an
element in $S_i$ with $h$-invariant conjugacy class.
 Finally the number of elements is a conjugacy class of $g$ in
 $S_i$ is $\vert C(g) \cap S_i \vert \over \vert S_i \vert$. Therefore,
 after splitting $X^{[S_i]}$ into a disjoint union of $X^{S_i}$ and omitting
terms corresponding to $S_i$ for which $h(S_i) \ne S_i$ (hence giving
zero contribution due to  the vanishing of the Lefschetz number),
 the latter sum can be replaced by
 the following sum over the pairs $(g,S_i),h(S_i)=S_i, g \in S_i $:
 $$\Sigma {{L(h,X^{[S_i]}) {\vert S_i \vert}
 \cdot \vert C(g) \cap S_i \vert}
 \over {\vert G \vert \cdot \vert S_i \vert }}=
 \Sigma_{(g,S_i),g \in S_i, h(S_i)=S_i}
  L(h,X^{S_i}) {{\vert C(g) \cap S_i \vert \cdot  \vert C(g) \vert}
\over {\vert C(g) \vert \cdot \vert G \vert }} \eqno (2.4)$$
Note that a centralizer $C(g)$ acts on the set of
stabilizers containing $g$. ($g \in S_i,c\in C(g)$ implies
 $g=cgc^{-1} \in cS_ic^{-1}$).
 The latter sum can be rewritten as a sum  over $C(g)$ orbits $C(g)(S_i)$
 of $S_i$:
 $$\Sigma L(h, \coprod X^{S_i}) {{\vert C(g) \cap S_i \vert} \over
 {\vert C(g) \vert}}
\cdot {{\vert C(g) \vert} \over {\vert G \vert}} \eqno (2.5)$$
 since the centralizers of $g$ in each subgroup in $C(g)(S_i)$ are
the same. This is equal to:
 $$\Sigma_{g,C(g)(S_i)}
  L(h,X^{C(g)(S_i)} /C(g)) {{\vert C(g) \vert} \over {\vert G \vert}}=
  \Sigma L(h,X^g/C(g)) {{\vert C(g) \vert} \over {\vert G \vert}} \eqno (2.6)$$
 where the summation is over all elements $g$ which have $h$ invariant
 conjugacy class in any $h$ invariant stabilizer to which it belongs.
 The last sum finally is  equal to
 $$\Sigma_{g \in {\cal C}(h)} L(h,X^g/C(g)) \eqno (2.7)$$
 where the summation is over
 all {\it conjugacy classes} in $G$ of elements with the specified
 property.
\bigskip
 \bigskip  {\bf 3.Applications to the calculation of
the Lefschetz number of the actions on mirror manifolds.}
\bigskip
\bigskip Here, using the formula from preceding section, we shall
  calculate the Lefschetz number of the
involution induced on the mirror of quintic $V$ in ${\bf P}^4$
considered in [COGP]:
 $$Q(\lambda x_0,x_1,x_2,x_3,x_4)=
 x_0^5+x_1^5+x_2^5+x_3^5+x_4^5-5 \lambda x_0x_1x_2x_3x_4=0
\eqno (3.1)$$
 when the
involution $h: (x_0,x_1,x_2,x_3,x_4) \rightarrow (x_1,x_0,
x_2,x_3,x_4)$
is induced by the interchange of first two coordinates.
\par Recall that the mirror is a crepant resolution of $\widetilde {V/G}$
where $G={\bf Z}^3_5$ is (the maximal) faithfully-acting quotient
of the group of automorphisms of $V$ given by
$$(x_0,x_1,x_2,x_3,x_4) \rightarrow (\omega_5^{\alpha_0}x_0,
\omega_5^{\alpha_1}x_1, \omega_5^{\alpha_2}x_2, \omega_5^{\alpha_3}x_3,
\omega_5^{\alpha_4}x_4) \eqno (3.2)$$
 where $\Sigma \alpha_i \equiv 0 mod 5$ and $\omega_5$ is a non-trivial root
of unity of degree 5.
First note that the condition $g \in {\cal C}(h)$ implies that
$\alpha_0=\alpha_1$ in (3.2).  Each 1-dimensional fixed point set
 is the intersection of $V$ with the subspace given by vanishing
of two coordinates and is fixed by elements (3.2) for which three exponents
$\alpha_i$ coincide. Since there are 3 possibility to choose $\alpha_i$
coinciding with chosen $\alpha_0=\alpha_1$ we obtain 12 elements
(different from identity)
in ${\cal C}(h)$ having 1-dimensional fixed point set. The Lefschetz
number of $h$ acting on the quotient of each of these fixed point sets
 is  equal to 2.
\par Similarly one can consider the   0-dimensional fixed point sets.
 Each is the fixed point set  of an element (3.1)
with two pairs of equal components ($\alpha_i=\alpha_j$ for two
pairs of indices $(i,j)$). Since $\alpha_0=\alpha_1$ we see that there
are  $12$ elements in ${\cal C}(h)$ having 0 dimensional fixed point set.
  Since the quotient in case of 0-dimensional fixed point sets
has euler characteristic 2 and $h$ acts trivially on it we see that the
contribution of $g \ne id$ in the expression (1.1)
 for $L(h,\widetilde {V/G})$ is
 $2 \times 12+ 2 \times 12=48$.
Let us consider the  contribution of the identity element.
 The cohomology of $V/G$ can be identified with
 the $G$-invariant
part of the cohomology of $V$. Hence the even dimensional cohomology of
$V/G$ in each dimension has rank 1. In dimension 3 the $G$-invariant part
of $H^3 (V,{\bf C})$ is generated by the residues of meromorphic forms
$$\omega_{\lambda}={{\Sigma (-1)^i x_idx_0 \wedge...\hat {dx_i}...
\wedge dx_4} \over {Q(\lambda,x_0,x_1,x_2,x_3,x_4)}},
{{d \omega_{\lambda}} \over {d \lambda}},
 {{d^2 \omega_{\lambda}} \over {d \lambda^2}},
{{d^3 \omega_{\lambda}} \over {d \lambda^3}} \eqno (3.3)$$
  via Griffiths theory (cf.[COGP], [M]).
 These forms  are clearly $h$-anti-invariant.
  Hence  $L(h,V^{id}/G)=8$. Therefore $L(h,\widetilde {V/G})=
8+48=56$.
 On the other hand $L(h,V)=e(V^h)$. The fixed point set  of $h$
 acting in $V$ consists
 of the point $(1,-1,0,0,0)$ and the points of the quintic in the
hyperplane $x_0=x_1$. i.e. a non singular quintic surface. Its euler
characteristic is $(3-(-10)) \times 5+ (-10)=55$ (since quintic
surface in ${\bf P}^3$ is a 5-fold cyclic
cover of projective plane branched over
a plane quintic having genus $6$). Hence we obtain:
 $L(h,V)=L(h,\widetilde {V/G})(=56)$.
\par Now let us consider the action on the quintic
 $(3.1)$ of the transformation:
$$h: (x_0,x_1,x_2,x_3,x_4) \rightarrow (x_1,x_0,x_2,x_4,x_3) \eqno (3.4)$$
Here ${\cal C}(h)$ consists of elements (3.2) for which $\alpha_0=\alpha_1,
\alpha_3=\alpha_4$. There are 4 such elements different from identity
and the Lefschetz number of $h$ acting on the quotient of their fixed point
 sets is equal to 2. The contribution of the identity element now is zero
since now the action of $h$ on forms $(3.3)$ is trivial. On the
other hand the fixed point set of (3.4) consists of the union of the line
$(x_0,-x_0,0,x_3,-x_3)$ and the plane quintic curve which is the intersection
of (3.1) and $x_0=x_1, x_3=x_4$ (which has genus 6). Hence $L(h,V)=
-L(h,\widetilde V/G)(=-8)$.
\par Finally let us consider example from [LT] in which $V$ is a
complete intersection on ${\bf P}^5$ given by equations:
 $$x_1^3+x_2^3+x_3^3=3\lambda x_4x_5x_6,$$ $$ x_4^3+x_5^3+x_6^3=3 \lambda
 x_1x_2x_3 \eqno (3.5)$$
with the group $G_{81}$ of order 81 acting on (3.5) as follows:
$$(x_1,x_2,x_3,x_4,x_5,x_6) \rightarrow
(\zeta_3^{\alpha_1} \zeta_9^\mu x_1, \zeta_3^{\alpha_2} \zeta_9^{\mu} x_2,
\zeta_9^\mu x_3,\zeta_3^{\alpha_4} \zeta_9^{-\mu} x_4,
\zeta_3^{\alpha_5} \zeta_9^{-\mu} x_5, \zeta_9^{-\mu} x_6) \eqno (3.6)$$
where $\alpha_i \in {\bf Z}_3, i=1,2,3,4,\mu \in {\bf Z}_3$ and
 $ \mu \equiv \alpha_1+\alpha_2 \equiv \alpha_4+\alpha_5 mod 3$.
Let us consider the following involution:
$$h: (x_1,x_2,x_3,x_4,x_5,x_6) \rightarrow (x_2,x_1,x_3,x_5,x_4,x_6)
\eqno (3.7)$$
The condition $g \in {\cal C}(h)$ implies that in (3.6) we have
$\alpha_1=\alpha_2,\alpha_4=\alpha_5$ i.e.  $g \in {\cal C}(h)$ must
have  the form: $$(x_1,x_2,x_3,x_4,x_5,x_6) \rightarrow
(\zeta_3^{\alpha}\zeta_9^{\mu} x_1,\zeta_3^{\alpha}\zeta_9^{\mu}x_2,
\zeta_9^{\mu}x_3,\zeta_3^{-\alpha}\zeta_9^{-\mu}x_4,
\zeta_3^{-\alpha}\zeta_9^{-\mu}x_5,\zeta_9^{-\mu}x_6) \eqno (3.8)$$
where $\alpha \in {\bf Z}_3,\mu \in {\bf Z}_9$ and $2\alpha=\mu mod 3$.
Hence we have 9 elements in ${\cal C}(h)$ out of which 8 non identity elements
have zero dimensional fixed point
set. The Lefschetz number of $h$ on the quotient of each of these
zero dimensional fixed point sets  by $G_{81}$ is equal to 2
For $g =id$ the contribution in (1.1) is the Lefschetz number
 $L(h,V^{id}/G_{81})$ which is equal to zero.
This  follows from explicit expression for the forms
representing $G_{81}$-invariant cohomology classes on complete intersection
 (3.5) (cf. [LT] (8) on p.32) as was done above in the case of quintic.
Hence $L(h,\widetilde {V/G_{81}})=0+2 \times 8=16$. On the other hand
the fixed point set of $h$ acting on (3.5) consists of the line:
 $(x,-x,0,y,-y,0)$ and the intersection of complete intersection (3.5)
 with the $x_1=x_2, x_4=x_5$  has the Euler characteristic $-18$.
Hence $L(h,V)=-L(h,\widetilde {V/G_{81}}(=-16)$.
\par One can wonder if there is a physical reason for these
simple equalities $L(h,V)=-sign (h) L(h,\widetilde {V/G})$
 between Lefschetz numbers of automorphisms of
a Calabi Yau manifold and its mirror which came out in these examples..
\bigskip \bigskip {\bf 4.Actions on resolutions of 2-dimensional
singularities}
\bigskip First we shall consider several explicit examples (cf.[Sl]).
Let $G={\bf Z}_n $ be a cyclic subgroup of the torus of $SL_2 ({\bf C})$
consisting of the matrices of the form:
 $$ \left (\matrix{ \omega_n^a & 0 \cr
                    0 & \omega_n^b \cr} \right) \eqno (4.1)$$
 where $a+b \equiv 0 mod n$. The matrix:
      $$\left (\matrix{ 0 & 1 \cr
                    -1 & 0 \cr} \right) \eqno (4.2)$$
acts on $G$ by conjugation  provided $n$ is even.
The number of invariant conjugacy classes
is $2$ while only one component of the exceptional set of the resolution
is $h$ invariant (cf. [Sl] p. 76). Hence the corresponding Lefschetz
number is $2$ and so  $G^h=L(E,h)$.
\par If $\bar h$ is given by $$\left (\matrix{ 0 & 1 \cr
                    1 & 0 \cr} \right) \eqno (4.3)$$
then the number of $\bar h$-invariant elements of $G$ is 1 (resp. 2)
if $n$ is even (resp. odd). The action on the resolution has 1 or 2
fixed point in respective cases.
\par Let us consider the case of binary dihedral group corresponding to the
Dynkin diagram $D_r$.
Recall that $D_r$ is a ${\bf Z}_2$ extension of ${\bf Z}_{2(r-2)}$ which
is the
subgroup of $SL_2({\bf C})$ generated by the cyclic group
$G$ given by matrices $(4.1)$
with $n=2(r-2)$ and (4.2). Then $D_{2r-2}$ contains $D_r$ as a normal
subgroup of index 2. Let $h$ be the nontrivial element of the
quotient and let $C$ be $h$-invariant subgroup of $G$
($Card C=2$).
$D_r$ has $Card C$ conjugacy classes having 1 element (i.e., the elements
of $C$), $Card G-Card C/2$ classes containing 2 elements (of set $G-C$)
and 2 conjugacy classes union of  which forms the coset of $G$ in $D_r$.
Among these $r+1$ conjugacy classes $r-1$, corresponding to those
in first two groups are, invariant under the action of the group
${\bf Z}_{4(r-2)}/{\bf Z}_{2(r-2)}$ (all cyclic groups are
groups of matrices as above) i.e., we obtain $r-1$ invariant conjugacy classes.
 The action of the non trivial element $\tilde h$
of the latter on the resolution of ${\bf C}^2/D_r$ fixes
the chain on $r-2$ rational curves (cf. [Sl] p.76 ).
Hence the Lefschetz number is equal to
$r-1$ i.e. the number of invariant conjugacy classes is $L(h)$.
\par In the case of binary tetrahedral group (i.e. the extension of
 quaternionic group by ${\bf Z}_3$) out of 7 conjugacy classes 3 are invariant.
In the case $D_4$ one has additional automorphisms of order 3. It leaves
invariant $2$ conjugacy classes and acts on the resolution leaving fixed
one component of the Dynkin diagram.
\par In general one can deduce the equivariant version of McKay
correspondence (**)
from the geometric description of McKay correspondence due to
Gonzalez-Springer and Verdier (cf. [GSV]). The number of $h$-invariant
conjugacy classes can be identified with the number of $h$-invariant
representations of each of binary polyhedral group $B$. Since the extension
of the bundle on the minimal resolution of ${\bf C}^2/B$ corresponding
to an $h$-invariant representation of $B$ will be $h$-invariant, its
first Chern
class will be invariant i.e. corresponding exceptional set will be invariant
and the conclusion follows since the traces on $H^0$ and $H^2$ of the
 resolution are obvious.
\bigskip \bigskip {\bf 5.Symmetries of resolution of abelian quotient
singularities}
\bigskip Let $X$ be a toric variety i.e. a torus ${\bf T}$ acts on $X$
 and action of ${\bf T}$ on one of the orbits (which we  shall
denote ${\cal T}$) is transitive.
 Let $h$ be a biregular automorphism of $X$ normalizing
 ${\bf T}$.  Then $h$ acts on the torus ${\bf T}$  as follows:
 $h(t)$ is the unique element of the torus which takes any point
 of  $v$ of ${\cal T}$ into $h^{-1} t h (v)$.
\bigskip Now let ${\bf T}$ be the maximal torus of
 $SL_n({\bf C})$ consisting
of diagonal matrices having traces equal to $1$.
 Let $H$ be an abelian subgroup of ${\bf T}$.
Let $h$ be an element of $GL_n ({\bf C})$ which normalizes ${\bf T}$
and $H$. $h$ can be viewed as an automorphism of both ${\bf T}$  and  $H$.
$h$ also acts on the quotient ${\bf C}^n /H$. The latter is a toric
variety in a natural way and the lattices of 1-parameter subgroups $M$ and $N$
of the dense tori of ${\bf C}^n$ and ${\bf C}^n/H$
are related by the following sequence:
$$0 \rightarrow M \rightarrow N \rightarrow H \rightarrow 0 \eqno (5.1)$$
Let $l$ be the order of $h$ acting of the lattice $M$ (or $N$).
The quotient of the normalizer of ${\bf T}$
in $GL_n ({\bf C})$
is the symmetric group $S_n$. Let $(l_1,..,l_k)$
be the sequence of lengths of cycles of the permutation (of characters
given by coordinates) defined by $h$. One has $l=l.c.m. (l_1,...,l_k)$.
\par  Let $\Delta$ be the unit simplex $\{ (x_1,..,x_n) \subset
 M \otimes {\bf R} \vert x_1+...x_n \le 1 \}$. The only face of this
simplex which does not belong to any coordinate plane will be called
the base of this unit simplex.
\bigskip We are going to construct standard triangulations of
certain simplices which will be used below.
Let $M={\bf Z}^{l_1} \oplus ... \oplus {\bf Z}^{l_p}
 \oplus {\bf Z}^{m_1} \oplus ... \oplus {\bf Z}^{m_q}$
and let $e^i_j$ ($i=1,..,p, j=1,...,l_i $ or $i=1,,,q,j=1,...,m_q$)
be the standard generators of each direct summands.
 Let the action of $g$ on $M$ be given by the cyclic permutation of the vectors
of the standard basis of each summand:
 $g(e^i_j)=e^i_{s_i(j)}$ where $s_i$ is the cyclic
permutation of the integers $1,..,l_i$ or $1,...,m_i$.
 Let $L$ be the $k+\Sigma_{i=1}^{i=q} m_i$-dimensional subspace
given by the equations $x^i_a=x^i_b$ ($i=1,...,k, 1 \le a,b \le l_i$).
The volume of the simplex $L \cap \Delta$ is ${1 \over {l_1 \cdot \cdot \cdot
 l_k}} \cdot {1 \over {dim \Delta \cap L !}}$. Let us consider the following
triangulation of $\Delta$ by simplices of the form $\Delta(r_1,...,r_k)=
(...,\hat X^i_j(r_i),...a_0,..,a_{k})$ where $\hat X^i_j(r_i)$ is the
collection of points all coordinates of which but one are zeroes and
nonzero coordinate is 1 corresponds a vector $e^i_j$ where $j \ne r_i$,
$a_0$ is the origin and $a_i$ ($1 \le i \le k$) are the vertices of
 of $L \cap \Delta$. The action of $g$ on this triangulation is given by
 $g( (...,\hat X^i_j(r_i),...a_0,..,a_{k}))=
(...,\hat X^i_j(s_i(r_i)),...a_0,..,a_{k})$. The pair
 $(\Delta, \Delta \cap L)$ will be called the standard pair
of type $(l_1,..,l_s \vert m_1,...,m_t)$.
\bigskip {\bf Definition} Let $\Sigma$  be a simplex with the
vertices of a lattice $M$ on which $g \in Aut M$ acts
simplicially, $T$ is an $g$-invariant triangulation of $\Sigma$ and
$\Sigma^g$ be the fixed point set of $g$. $T$ is called standard if
the pair $(\Sigma,\Sigma^g)$ is isomorphic to $(\Delta,\Delta \cap L)$
by an isomorphism preserving volume of simplices of dimension
 $dim \Sigma^g+1$.
\bigskip  {\bf Remarks}
1. The vertices of simplex $\Sigma^g$ are not necessarily in $M$.
\par 2. The volume of each simplex is calculated according to the measure
induced by the lattice in the linear subspace supporting this simplex.
\bigskip The properties of $g$-standard simplices which will be used later
are the following:
\par a) $g$ acts on the vertices of a $g$-standard simplices by
permutation with cycles of length $(l_1,...,l_s)$.
\par b) Codimension of $\Sigma \cap L$ in $\Sigma$ is $\Sigma (m_i-1)$.
\par c) $g$ acts on simplices of dimension $dim \Sigma \cap L+1$ as
permutation having the lengths of cycles $(m_1,..m_t)$.
\par d) the volume of $\Sigma \cap L$ is
${1 \over {l_1 \cdot \cdot \cdot l_s}} \cdot {1 \over {dim \Sigma^g !}} $.
\bigskip {\bf Definition.} A triangulation of a simplex is $g$-adjusted if
it is a refinement of a triangulation in which each $g$-invariant simplex
is $g$-standard.
\bigskip {\bf Theorem 2.} If there exist an $h$-adjusted triangulation
 of the unit simplex
$\Delta$ in $M$  such that each simplex of triangulation has vertices
in $N \cap \Delta$ and such that the volume of each simplex
relative to the lattice induced by $N$ on the linear subspace supporting
its simplex is 1 then there exist  a $h$-invariant crepant resolution of
 ${\bf C}^n/H$ such that the Lefschetz number of $h$ acting on the
latter is equal to the order of the group of $h$-invariant elements of $H$.
\bigskip {\bf Proof.} Since any triangulation of $b\Delta$ with vertices in
$N$ induces a crepant resolution (cf. [R],[BD]), we obtain a crepant
resolution from an $h$-adjusted triangulation
mentioned in the statement.  Let $L(h)$ be  the Lefschetz number of
$h$ acting on a chosen resolution, $Card H^h$ be the number of
 $h$-invariant elements of $H$ and  $(l_1,...,l_k)$ be the sequence
of greater than $1$ lengths of cycles of permutation corresponding to $h$.
 Let $L$ be the subspace of $M \otimes {\bf R}$ of elements fixed
by $h$. We claim that one has
\par 1. $L(h)=l_1 \cdot \cdot \cdot l_k \cdot vol_{N \cap b\Delta}
  (L \cap b\Delta) \cdot k!$
\par 2. $vol_{N \cap b\Delta}={1 \over {l_1 \cdot \cdot \cdot l_k}}
 \cdot {1 \over k!}  Card H^h$.
\par Clearly the theorem follows from 1 and 2.
\par Let us first calculate the Lefschetz number of $h$ acting on the
resolution. $h$ acts freely outside of the union of tori
corresponding to the simplices which intersect $L$, since any simplex
fixed by $h$ contains a fixed point of $h$ and  $L$ is the total set of fixed
points of $h$. Suppose that a simplex $\sigma$ is $h$ invariant
and not in the closure of an invariant simplex.
Then since triangulation is $h$-adjusted, its vertices
 permuted by a permutation consisting of cycles
of lengths  $(l_{i_1},...,l_{i_k})$  and the union ${\cal B}_ {\sigma}$
of certain simplices
 represent a collection of simplices of dimension $dim \sigma+1$
permuted by permutation with cycles with sequence  of lengths
complementary to $( l_{i_1},...,l_{i_k})$. Let ${\cal A}_{\sigma}$
be the union of these simplices and let $A_{\sigma}$ and $B_{\sigma}$
be the corresponding toric varieties. The Lefschetz number of
$h$ acting on both $A_{\sigma}$ and $B_{\sigma}$ are equal to the Lefschetz
number of $h$ acting on the torus corresponding to $\sigma$. On the
other hand viewing the torus ${\bf T}_{\sigma}$
as a subset of $B_{\sigma}$ shows that the matrix of $h$ acting on the
$H_1 ({\bf T}_{\sigma},{\bf Z})$ is formed by the  blocks
${\bf A}_{l_1-1},...,{\bf A}_{l_k-1}$ where ${\bf A}_s$ is the
following $s \times s$ matrix:
$$  \left (\matrix{ 0 & 1 & 0 & ...& 0 \cr
                    0 & 0 & 1 & ... &0 \cr
                    . & . & . & .   & 0 \cr
                    -1& -1 &...& ...      &- 1 \cr} \right) $$
If $s > 1$ then $det(I-A_s)=s+1$ (and $0$ for $s=1$).
If $\phi$ is transformation  of a torus ${\bf T}$ and $\phi_*$ is the
corresponding automorphism on $H_1 ({\bf T},{\bf Z})$ then the Lefschetz
number of $\phi$ is equal to $det (I-\phi_*)$ (since
 $H^*({\bf T},{\bf Z})=\Lambda ^*(H_1({\bf T},{\bf Z})$ and the
eigenvalues of the map induced by $\phi$ on $\Lambda^k H_1({\bf T},{\bf Z})$
 are the elementary symmetric functions of the eigenvalues of
$\phi$ acting on $H_1 ({\bf T},{\bf Z})$). Hence the
Lefschetz number of $h$ acting on ${\bf T}_{\sigma}$ is $l_1 \cdot \cdot
\cdot l_k$. This calculation implies that the total Lefschetz number
of $h$ is $$\Sigma {{l_1 \cdot \cdot \cdot l_k} \over
{l_{i_1} \cdot \cdot \cdot
 l_{i_s}}} r_{l_{i_1},...,l_{i_s}}$$
 where $ r_{l_{i_1},...,l_{i_k}} $ is the total number of
$h$ adjusted simplices of type $(l_{i_1},...,l_{i_s})$.
 \par On the other hand the volume of $L \cap b\Delta$
  is equal to $$ 1/ k! \Sigma {1 \over {l_{i_1} \cdot \cdot \cdot
 l_{i_s}}} r_{l_{i_1},...,l_{i_k}}$$
 since for a  $h$ standard simplex $\Sigma$ of type $(l_{i_1},..,l_{i_s}$
 for which $dim \Sigma \cap L=k$ is equal to $1 \over {l_{i_1} \cdot
\cdot \cdot l_{i_k}}$.
 Hence the claim 1 above follows.
\par The number $Card H^h$ is the volume of the unit parallelepiped of
lattice $M \cap L$ measured relative to the lattice $N \cap L$. It is
clear from direct calculation that $vol L \cap b\Delta =1/l_1...l_k
 vol b\Delta_M$. This shows 2.
\bigskip {\bf Remark.}
 Let us consider an orbit corresponding to a simplex $\sigma$
in the closure of an invariant simplex $\tau$. We claim that
the Lefschetz number of $h$ acting on the orbit corresponding
 to $\sigma$ is equal to zero. We shall show that the orbit
 corresponding to $\sigma$  contains one dimensional orbit fixed
 by $h$ pointwise.
  We can assume that $dim \sigma= dim \tau -1$. Let $X_{\sigma}$
 be the closure in $X$ of the orbit corresponding to $\sigma$. Then
$h$ acts on $X_{\sigma}$ and the latter contains $h$ invariant codimension
one orbit (the orbit corresponding to $\tau$). Let us consider the
fan corresponding to this toric variety and an orbit $\phi$ corresponding
to the point of the lattice which belong to $h$-invariant ray corresponding
to $h$ invariant codimension one orbit. This ray is fixed by $h$ and hence
the lattice point on it is fixed by $h$. We claim that this one
dimensional orbit is fixed pointwise. Indeed  $h$ is a transformation
which fixes limit of  this orbit when $t \rightarrow 0$
(note that $h$ is either the identity or sends $x \rightarrow x^{-1}$).
\bigskip {\bf  Theorem 3.} An $h$-adjusted triangulation of the unit simplex
 exist in dimension 3. In particular the conclusion of theorem
2 is always  true in dimension 3.
\bigskip {\bf Proof.} The order of a nontrivial automorphism
$h \in S_3$ is either 2 or 3. If $ord h=2$, $h$ fixes vertex $P$
and the intersection of $L$ with the side $O_P$ of $b \Delta$ opposite to
$P$ belongs to the lattice $N$, then we can split $b \Delta$ by $L$
into the union of two triangles,  $T$ and $hT$,
 take  a triangulation of $T$
in which all vertices are in $N \cap b \Delta $ and have area $1$ and then
take $h$ image of this triangulation in $hT$. If $ord h=2$ but
 $L \cap O_P$ is not in $N$,
 consider the triangle $T_1$ of area one with vertex
in the point $L \cap N$ closest to $O_P$ and two vertices of $N$ on $O_P$,
then triangulate $T - T  \cap T_1$ by triangles of area $1$ with vertices
in $N$ and take its $h$ image to triangulate $hT-hT \cap T_1$. Now any segment
along $L \cap b \Delta$ is $h$-standard and $T_1$ is $h$-standard of type
$(2)$. Hence the triangulation is $h$ -adjusted.
\par If $ord h=3$ then one takes triangulation of one of the triangles
 formed by $L \cap b \Delta$ one then extends it to triangulation
of $b \Delta$ using $h$.
\bigskip {\bf Example.} Let us consider the action of subgroup
 ${\bf Z}_5^2$ of $\bf T$
consisting of matrices of the form:
  $$  \left (\matrix{ \omega_5^a &  0 & 0 \cr
                    0 & \omega_5^b & 0 \cr
                    0 & 0  &\omega_5^c \cr} \right) $$
where $a+b+c \equiv 0 mod 5$  where $\omega_5$ is a root of unity of
 degree $5$. $h \in GL_3 ({\bf C})$ given by
 $$  \left (\matrix{ 0 &  1 & 0 \cr
                    0 & 0 & 1 \cr
                    1 & 0  & 0 \cr} \right) $$
induces cyclic permutation of the entries of the matrices of ${\bf T}$.
The Lefschetz number of the automorphism of the resolution of
 ${\bf C}^3/{\bf Z}^2$ constructed above is equal to $1$ (coming from
the only invariant simplex which the cone over triangle having $(1/3,1/3,1/3)$
in its center and for which the length of intersection with $L$ is $1/3$)
and the only $h$ invariant element of ${\bf Z}_5^2$ is the identity.
\bigskip
\centerline {\bf References}
\bigskip
\bigskip [BD] {\bf V.Batyrev, D.Dais}, Strong McKay correspondence, String
theoretical Hodge Numbers and Mirror Symmetry. Preprint, 1994.
\bigskip [Br] {\bf R.Brown}, Lefschetz fixed point theorem,
 Scott, Foresman and Co,
Glenview, Ill, London. 1971.
\bigskip [COGP] {\bf P.Candelas, X.de la Ossa, P.Green, L.Parkes}, A pair
of Calabi Yau manifolds as an exactly soluble superconformal theory,
Nucl. Phys. ${\bf B}$ 359, 1991. 21-74.
\bigskip [GSV] {\bf G.Gonzalez-Springer and J.L.Verdier},
Construction Geometric
de la correspondence de McKay, Ann. Sci. Ecole Normale Sup. 16 (1983),
409-449.
\bigskip [LT] {\bf A.Libgober, J.Teitelbaum}, Lines on Calabi Yau complete
intersections, mirror symmetry and Picard Fuchs equations, Intern. Math.
Research Notices, No.1 1993, p.29-39.
\bigskip [M] {\bf D.Morrison}, Mirror Symmetry and rational curves on quintic
threefolds: a guide for mathematicians. Journ. of Amer. Math. Soc. 1993.
vol. 6,223-247.
\bigskip [R] {\bf S.Roan}, On the generalization of Kummer surfaces,
J. of  Differential Geometry, 30 (1989), 523-537.
\bigskip [Sl] {\bf P.Slodowy} Simple Singularities and Simple Algebraic
Groups. Lecture Notes in Math.  vol. 815. Springer Verlag, 1980.
\end